\documentclass[12pt]{amsart}

\theoremstyle{definition}

\newtheorem{theorem}{Theorem}[section]

\theoremstyle{definition}

\theoremstyle{remark}

\numberwithin{equation}{section}

\newcommand\Calderon{Cal\-der{\'o}n}

\newcommand{\cinf}[1]{C^\infty(#1)}

\newcommand{\dom}{\operatorname{dom}}
\newcommand{\End}{\operatorname{End}}

\newcommand{\ind}{\operatorname{ind}}

\newcommand{\rrr}{\!\!\upharpoonright\!} 
\newcommand{\sign}{\operatorname{sign}}
\newcommand{\wt}{\widetilde}

\def\ftoo#1{\stackrel{#1}{\too}}

\newcommand{\dd}[1]{\frac{d}{d#1}}

\def\gS{\Sigma}

\def\lla{\langle}
\def\noi{\noindent}

\def\rra{\rangle}

\def\too{\longrightarrow}

\def\wt{\widetilde}

\def\={\cong}
\def\>{\supset}
\def\<{\subset}
\def\ii{^{-1}}
\def\12{\frac{1}{2}}
\def\0{^{\circ}}

\usepackage[english]{babel}

\usepackage[latin1]{inputenc}

\usepackage{times}
\usepackage[T1]{fontenc}

\usepackage{graphicx}

\usepackage{amscd}
\usepackage{amssymb}

\begin{document}

\noi {\bf Continuous Variation of Cauchy Data Spaces}.
{\bf Task - Results - Open Problems}

\smallskip


\noi 
B. Boo{\ss}-Bavnbek, Roskilde University, Denmark\\
Matthias Lesch, University Bonn, Germany, \\
Chaofeng Zhu, Chern Institute, Tianjin, China

\smallskip

\noi International Seminar {\em Analysis of spectral invariants and related operator theory,
Tokyo University of Science, Unga Campus, 5-6 Oct. 2009}
\bigskip

Over the last decades, substantial progress has been achieved 
investigating spectral properties of geometrically defined differential operators
like the Laplacian and the Dirac operator on manifolds with boundary and on partitioned manifolds.

Some common deep functional-analytical roots of these formulas have been revealed in
[B. Boo{\ss}-Bavnbek and K. Furutani, Tokyo J. Math. 21 (1998), no. 1, 1--34; MR1630119 (99e:58172)],
emphasizing the role of the Cauchy data spaces. More precisely, 
in the von Neumann-Kre\u{i}n-Vishik-Birman tradition one is given
a complex separable Hilbert space $\mathcal{H}$ and  a closed symmetric operator $A$. One
defines the symplectic Hilbert space of abstract boundary values by 
$\beta(A):=\dom(A^*)/\dom(A)$ with naturally induced inner product 
$\lla [x],[y]\rra $ and symplectic form $\omega([x],[y])$ and the 
natural {Cauchy data space} $\operatorname{CD}(A):= \{[x]\mid x\in\ker A^*\}$. 
One has a canonical correspondence between all self-adjoint extensions $A_D$ 
of $A$ with domain $D$ and the Lagrangian subspaces $[D]\< \beta(A)$. In this framework, e.g.,
if $A_D$ is
a {self-adjoint Fredholm extension} and  $\{C_t\}$ a {continuous curve} in $\mathcal{B}(\mathcal
{H})$ with $\ker(A^*+C_t+s)\cap \dom(A)=\{0\}$ for small $|s|$ ({weak inner
UCP}), one obtains that $\{\operatorname{CD}(A+C_t),[D]\}$ is a continuous curve of Fredholm pairs
of Lagrangians and $\operatorname{SF}\{(A+C_t)_D\} =
\operatorname{MAS}\{\operatorname{CD}(A+C_t),[D]$, relating the spectral flow of a self-adjoint
Fredholm operator under bounded variation with the Maslov index of the corresponding curve of
Lagrangians in the abstract boundary space.

The strength of this functional-analytical approach shows up when dealing
with systems of ordinary differential equations on the
  interval, generalizing the classical Morse index theorem for geodesics on Riemannian manifolds 
to Subriemannian manifolds. It recovers  the Floer-Yoshida-Nicolaescu splitting results for the
spectral flow of curves of Dirac operators on partitioned manifolds
(i.e., the family version of the Bojarski Conjecture), and it 
provides a basic functional-analytical model for {\em
quantization} and {\em tunneling}, relating {\em spectral} and {\em
symplectic} invariants. 

In our joint papers, we wanted to transgress the limitations of the general functional-analytical
approach: we don't wish to keep the domain  fixed under variation; nor to restrict to bounded 
(i.e., 0 order) perturbations; nor to confine the applicability to ordinary differential
 equations or Dirac   type operators with constant coefficients in normal direction
  (product case) close to boundary. We took up the challenge of 
general linear elliptic differential operators; investigated weak inner UCP;
established the existence of self-adjoint Fredholm extensions; 
admitted variation of domain and skew boundaries, and investigated
uniform structures and continuous perturbations.
We conclude that the ``natural" (von Neumann) approach is insufficient, and
more {analysis} (splitting the coefficients near the boundary and 
{pseudodifferential} calculus) are neded.

Let $M$ be a smooth compact Riemannian manifold with boundary $\gS$,
$E,F$ Hermitian vector bundles over $M$, and
$A:C^\infty(M,E)\to C^\infty(M,F)$ an elliptic differential
operator (of first order). Recall that 
$\rho:L^2_s(M,E)\to L^2_{s-1/2}(\gS,E_{\gS})$ for $s> 1/2$ is
extendable to ${\mathcal D}_{\operatorname{max}}(A)$. Then the classical definition 
of the {Cauchy data space} $N_+^0(A)$ of $A$ is the closure of$\{\rho u\mid Au=0 
\text{ in } M\setminus\gS,\, u\in C^\infty(M,E)\}$ in ${L^2(\gS,E_{\gS})}$.  
In [Amer. J. Math. 88 (1966) 781--809; MR0209915 (35 {\#}810) and
{\em Pseudo-Diff. Operators (C.I.M.E., Stresa, 1968)}, pp. 167--305 Edizioni Cremonese, Rome 1969; 
MR0259335 (41 {\#}3973)], R.T. Seeley proved that this Cauchy data space can be obtained as the range
of a pseudodifferential projection. The basic ingredients for Seeley's result have 
been the construction of an invertible extension $\wt A$ of $A$ over a closed manifold $\wt M$
by extending $A$ to a collar, then doubling and applying symbolic calculus and UCP management.
As a result, he received a Poisson operator $K_\pm:=\pm r^\pm \wt A\ii\rho^* J(0)$ where 
$J(0)=\sigma(A)(\cdot,\nu)\in \End(E_{\gS})$ denotes the principal symbol of $A$
in normal direction at the boundary. He showed that the operator $C_\pm := \rho K_\pm$ 
is a pseudodifferential projection onto $N^0_+(A)$ and called it the {Calder{\'o}n projection}.
 
In [{\em Elliptic boundary problems for Dirac operators}, 
Birkhäuser, Boston, 1993; MR1233386 (94h:58168)],
it was shown by B. Boo{\ss}-Bavnbek and K.P. Wojciechowski that Seeley's construction is canonical for Dirac type operators in product form close to the boundary 
and yields the Lagrangian property of the Cauchy data space. The reason is that for 
such operators the invertible extension $\wt A$ can be explicitly defined on the true closed double 
$\wt M$ of $M$ and does not involve any choices. As a consequence, the Cauchy data spaces, respectively, the
\Calderon\ projection vary continuously under smooth deformation of the data defining the
Dirac operator [B. Boo{\ss}-Bavnbek, M. Lesch and J. Phillips, Canad. J. Math. 57 (2005), no. 2, 225--250; 
MR2124916 (2006a:58029)].

We traced these results back to, what we call 
``{Dirac Operator Folklore}": (i) weak inner UCP, i.e.,
$\ker A\cap {\dom}(A_{\operatorname{min}})=\{0\}$ with
${\dom}(A_{\operatorname{min}})=L^2_{1,\operatorname{comp}}(M,E)$; 
(ii) symmetric principal symbol of the tangential operator $B$ in the
decomposition $A=J_0(\partial_x+B)$ where $x$ denotes the inner normal variable; and
(iii) a precise invertible double. From that alone, we emphasize, one can derive the 
transparent definition of the {Calder{\'o}n projection}, the Lagrangian property
of the Cauchy data space, the existence of a self-adjoint Fredholm extension given by 
a regular pseudodiffernetial boundary condition, the Cobordism Theorem, and the 
continuous dependence of input data.

We ask 
{\em How special are operators of Dirac type compared to arbitrary linear 
first order elliptic differential operators?} To answer that question, 
we first bring a given elliptic differential operator of first order in product form
$A = {J\bigl(\partial_{x}+B\bigr)}$
close to the boundary by suitable choice of the metric. 
Here, $J$ and $B$ vary with the normal variable $x$.
Note that dropping the geometric 
Dirac operator context, the metric structures need no longer to be fixed. 

Next, we give a canonical construction of a
{new invertible double} $\tilde A_{T}$
with $\mathcal D(\tilde A_{T}):=\{\binom{e}{f}\in L^2_1(M,E\oplus
F)\mid \varrho f=T\varrho e\}$ where
$\tilde{A} : C^\infty(M,E\oplus F)\ftoo{A\oplus(-A^t)}
C^\infty(M,F\oplus E)$ and 
$T \in\operatorname{Hom}(\gS,E|_\gS,F|_\gS)$ invertible
bundle homomorphism with  $J_0^*T$ positive definite. Then
$\tilde{A}_{T}$ is a {Fredholm operator with compact resolvent} with 
$\ker\tilde{A}_{T} =
Z_{+,0} \oplus Z_{-,0}$ and $\operatorname{coker}
\tilde{A}_{T} \simeq Z_{-,0} \oplus Z_{+,0}$ where
$Z_{+,0} := \{f\in L^2_1(M,E)\mid Af=0, \varrho f=0\}$ and
$Z_{-,0}$ denotes the corresponding kernel of $A^t$\,. 
For the most part of our work we pick 
$T := (J_0^t)^{-1}$\,. Denoting the pseudo--inverse of
$\tilde{A}_{T}$ by  $\tilde{G}$ , we define Poisson operators
$K_{\pm}  :=  \pm\, r^{\pm}\tilde{G}\varrho^*J_0:
L^2_s(\Sigma,E)\to L^2_{s+\12}(M,E)$ $(L^2_{s+\12}(M,F))$
and \Calderon\ operators $C_{+} := \varrho_+K_+, \quad C_{-}  :=
T^{-1}\varrho_{-}K_{-}$. We obtain that 
$C_{\pm}$ are projections with $C_+ + C_{-} = I$ and
$C_+(L^2)  =  N_+^0,\  C_-(L^2)   = T^{-1} N_-^0$\,.

The most delicate part of our work has been the investigation of the mapping properties of
the pseudo-inverse $\wt G$, the Poisson operators $K_\pm$ and the
\Calderon\ projection $C_\pm$\,.

Our model operator is
$A=J\bigl(\frac{\partial }{\partial x} +B(x))+0. \text{order}$. From the
ellipticity of $A$ we have that $i\xi+B(x)$ is invertible for real $\xi$
of sufficiently large numerical value (array of minimal growth). We put
$Q_+(x):=\frac{1}{2\pi i} \int_{\Gamma_+} e^{-x\lambda} (\lambda- B_0)^{-1}d\lambda$
a family of sectorial projections where $\Gamma_+$ is a contour which encircles the eigenvalues of $B_0$ in the right half plane. We
notice that $Q_+(x)$ corresponds $e^{-x B_0}1_{[0,\infty)}(B_0)$ if $B_0=B_0*$.
We had to display a delicate balance on a knife edge between 
general operator theory and pseudodifferential calculus when we realized that 
a priori $Q_+(x)=O(\log x), x\to 0+$, hence $P_+:=Q_+(0)$ is possibly unbounded. Within the pseudodifferential calculus, it follows, however, from
[T. Burak, Ann. Scuola Norm. Sup. Pisa (3) 24 1970 209--230; MR0279633 (43 \# 5354),
K. Wojciechowski, Simon Stevin 59 (1985), no. 1, 59--91; MR0795272  (86k:58120),
R. Ponge, J. Reine Angew. Math. 614 (2008), 117--151; MR2376284  (2008j:58039)]
that $P_+:=Q_+(0)$ is a bounded pseudodifferential projection.  
A postiori, we obtain $Q_+(x)\to P_+$ strong, $x\to 0+$.

Another hopefully useful concept introduced by us is
{the approximative Poisson operator}
$R:C^\infty(\partial M,E)\longrightarrow C^\infty(\mathbb{R}_+\times
\partial M,E\oplus F)$ with
$R\xi (x):=\varphi(x){Q_+(x)\xi\choose T Q_-(x)\xi}$, where
$\varphi$ is a suitable { cut-off function at } 0. We find
$R=\tilde A_T^{-1}\varrho^*$ +{ regularising remainder } and analyze the mapping property
of $R:L_s^2(\partial M,E)\to L^2_{s'}(\mathbb{R}_+\times \partial
M,E\oplus F)$ in dependence of $A$.

Regarding uniform structures, we find that
$C_+(A)- P_+(B_0)$ is a pseudodifferential operator of order $-1$ and that
$A\mapsto C_+(A)$ is as regular as $A\mapsto P_+(B(0))$
under the condition $\dim Z_0(A), \dim Z_0(A^t)=\text{const}$.
Moreover we show that $(A,P)\mapsto A_P$ is continuous in graph topology, if $P$ runs
in the space of ``regular" boundary conditions.

Further applications for $A=A^t$ are that the Cauchy data space is Lagrangian in 
the hermitian symplectic
Hilbert space $L^2(\partial M,E), \langle\cdot,J(0)\cdot\rangle)$; the
existence of a self-adjoint Fredholm extension $A_{C_\pm}$ (for suitable choice of $T$);
and the  cobordismus invariance of the index for arbitrary symmetric elliptic differential operators
on closed manifolds: $\sign iJ(0)$ vanishes on
$\bigoplus\limits_{\lambda \text{ imaginary}}\ker
(B(0)-\lambda)^N$, $N>>0$. 

The results are inspired by [B. Himpel, P. Kirk and M. Lesch, Proc. London Math. Soc. (3) 89 (2004), no. 1, 241--272; MR2063666 (2005f:58032)] and have been announced in
[B.Boo{\ss}-Bavnbek and M. Lesch, Lett. Math. Phys. 87 (2009), no. 1-2, 19--46; MR2480643].
The details are worked out in [B. Boo{\ss}-Bavnbek, M. Lesch and C. Zhu, 
J. Geom. Phys. 59 (2009), no. 7, 784--826].

\end{document}